\documentclass[12pt]{article}

\catcode`\@=11

\thicklines
\newskip\Einheit \Einheit=.6cm
\newcount\xcoord \newcount\ycoord
\newdimen\xdim \newdimen\ydim \newdimen\PfadD@cke \newdimen\Pfadd@cke
\def\PfadDicke#1{\PfadD@cke#1 \divide\PfadD@cke by2 
\Pfadd@cke\PfadD@cke \multiply\PfadD@cke by2}
\long\def\LOOP#1\REPEAT{\def\BODY{#1}\ITERATE}
\def\ITERATE{\BODY \let\next\ITERATE \else\let\next\relax\fi \next}
\let\REPEAT=\fi
\def\Punkt{\hbox{\raise-2pt\hbox to0pt{\hss\scriptsize$\bullet$\hss}}}

\def\DuennPunkt(#1,#2){\unskip
  \raise#2 \Einheit\hbox to0pt{\hskip#1 \Einheit
          \raise-1.5pt\hbox to0pt{\hss\tiny$\bullet$\hss}\hss}}
		  
\def\NormalPunkt(#1,#2){\unskip
  \raise#2 \Einheit\hbox to0pt{\hskip#1 \Einheit
          \raise-3pt\hbox to0pt{\hss\large$\bullet$\hss}\hss}}
\def\DickPunkt(#1,#2){\unskip
  \raise#2 \Einheit\hbox to0pt{\hskip#1 \Einheit
          \raise-4pt\hbox to0pt{\hss\Large$\bullet$\hss}\hss}}
\def\Kreis(#1,#2){\unskip
  \raise#2 \Einheit\hbox to0pt{\hskip#1 \Einheit
          \raise-4pt\hbox to0pt{\hss\Large$\circ$\hss}\hss}}
\def\Diagonale(#1,#2)#3{\unskip\leavevmode
  \xcoord#1\relax \ycoord#2\relax
      \raise\ycoord \Einheit\hbox to0pt{\hskip\xcoord \Einheit
         \unitlength\Einheit
         \line(1,1){#3}\hss}}
\def\AntiDiagonale(#1,#2)#3{\unskip\leavevmode
  \xcoord#1\relax \ycoord#2\relax \advance\xcoord by -0.05\relax
      \raise\ycoord \Einheit\hbox to0pt{\hskip\xcoord \Einheit
         \unitlength\Einheit
         \line(1,-1){#3}\hss}}
\def\Pfad(#1,#2),#3\endPfad{\unskip\leavevmode
  \xcoord#1 \ycoord#2 \thicklines\ZeichnePfad#3\endPfad\thinlines}
\def\ZeichnePfad#1{\ifx#1\endPfad\let\next\relax
  \else\let\next\ZeichnePfad
    \ifnum#1=1
      \raise\ycoord \Einheit\hbox to0pt{\hskip\xcoord \Einheit
         \vrule height\Pfadd@cke width1 \Einheit depth\Pfadd@cke\hss}%
      \advance\xcoord by 1
     \else\ifnum#1=2
      \raise\ycoord \Einheit\hbox to0pt{\hskip\xcoord \Einheit
         \unitlength\Einheit
         \line(0,1){1}\hss}
      \advance\xcoord by 0
      \advance\ycoord by 1
 \else\ifnum#1=3
      \raise\ycoord \Einheit\hbox to0pt{\hskip\xcoord \Einheit
         \unitlength\Einheit
         \line(1,1){1}\hss}
      \advance\xcoord by 1
      \advance\ycoord by 1
    \else\ifnum#1=4
      \raise\ycoord \Einheit\hbox to0pt{\hskip\xcoord \Einheit
         \unitlength\Einheit
         \line(1,-1){1}\hss}
      \advance\xcoord by 1
      \advance\ycoord by -1
   \else\ifnum#1=5
      \raise\ycoord \Einheit\hbox to0pt{\hskip\xcoord \Einheit
         \unitlength\Einheit
         \line(2,1){2}\hss}
      \advance\xcoord by 2
      \advance\ycoord by 1
	  \else\ifnum#1=6
      \raise\ycoord \Einheit\hbox to0pt{\hskip\xcoord \Einheit
         \unitlength\Einheit
         \line(2,-1){2}\hss}
      \advance\xcoord by 2
      \advance\ycoord by -1
	  \else\ifnum#1=7
      \raise\ycoord \Einheit\hbox to0pt{\hskip\xcoord \Einheit
         \unitlength\Einheit
         \line(3,1){3}\hss}
      \advance\xcoord by 3
      \advance\ycoord by 1
	  \else\ifnum#1=8
      \raise\ycoord \Einheit\hbox to0pt{\hskip\xcoord \Einheit
         \unitlength\Einheit
         \line(3,-1){3}\hss}
      \advance\xcoord by 3
      \advance\ycoord by -1
    \fi\fi\fi\fi\fi\fi\fi\fi
  \fi\next}
\def\hSSchritt{\leavevmode\raise-.4pt\hbox 
to0pt{\hss.\hss}\hskip.2\Einheit
  \raise-.4pt\hbox to0pt{\hss.\hss}\hskip.2\Einheit
  \raise-.4pt\hbox to0pt{\hss.\hss}\hskip.2\Einheit
  \raise-.4pt\hbox to0pt{\hss.\hss}\hskip.2\Einheit
  \raise-.4pt\hbox to0pt{\hss.\hss}\hskip.2\Einheit}
\def\vSSchritt{\vbox{\baselineskip.2\Einheit\lineskiplimit0pt
\hbox{.}\hbox{.}\hbox{.}\hbox{.}\hbox{.}}}
\def\DSSchritt{\leavevmode\raise-.4pt\hbox to0pt{%
  \hbox to0pt{\hss.\hss}\hskip.2\Einheit
  \raise.2\Einheit\hbox to0pt{\hss.\hss}\hskip.2\Einheit
  \raise.4\Einheit\hbox to0pt{\hss.\hss}\hskip.2\Einheit
  \raise.6\Einheit\hbox to0pt{\hss.\hss}\hskip.2\Einheit
  \raise.8\Einheit\hbox to0pt{\hss.\hss}\hss}}
\def\dSSchritt{\leavevmode\raise-.4pt\hbox to0pt{%
  \hbox to0pt{\hss.\hss}\hskip.2\Einheit
  \raise-.2\Einheit\hbox to0pt{\hss.\hss}\hskip.2\Einheit
  \raise-.4\Einheit\hbox to0pt{\hss.\hss}\hskip.2\Einheit
  \raise-.6\Einheit\hbox to0pt{\hss.\hss}\hskip.2\Einheit
  \raise-.8\Einheit\hbox to0pt{\hss.\hss}\hss}}
\def\SPfad(#1,#2),#3\endSPfad{\unskip\leavevmode
  \xcoord#1 \ycoord#2 \ZeichneSPfad#3\endSPfad}
\def\ZeichneSPfad#1{\ifx#1\endSPfad\let\next\relax
  \else\let\next\ZeichneSPfad
    \ifnum#1=1
      \raise\ycoord \Einheit\hbox to0pt{\hskip\xcoord \Einheit
         \hSSchritt\hss}%
      \advance\xcoord by 1
    \else\ifnum#1=2
      \raise\ycoord \Einheit\hbox to0pt{\hskip\xcoord \Einheit
        \hbox{\hskip-2pt \vSSchritt}\hss}%
      \advance\ycoord by 1
    \else\ifnum#1=3
      \raise\ycoord \Einheit\hbox to0pt{\hskip\xcoord \Einheit
         \DSSchritt\hss}
      \advance\xcoord by 1
      \advance\ycoord by 1
    \else\ifnum#1=4
      \raise\ycoord \Einheit\hbox to0pt{\hskip\xcoord \Einheit
         \dSSchritt\hss}
      \advance\xcoord by 1
      \advance\ycoord by -1
    \fi\fi\fi\fi
  \fi\next}
\def\Koordinatenachsen(#1,#2){\unskip
 \hbox to0pt{\hskip-.5pt\vrule height#2 \Einheit width.5pt depth1 
\Einheit}%
 \hbox to0pt{\hskip-1 \Einheit \xcoord#1 \advance\xcoord by1
    \vrule height0.25pt width\xcoord \Einheit depth0.25pt\hss}}
\def\Koordinatenachsen(#1,#2)(#3,#4){\unskip
 \hbox to0pt{\hskip-.5pt \ycoord-#4 \advance\ycoord by1
    \vrule height#2 \Einheit width.5pt depth\ycoord \Einheit}%
 \hbox to0pt{\hskip-1 \Einheit \hskip#3\Einheit 
    \xcoord#1 \advance\xcoord by1 \advance\xcoord by-#3 
    \vrule height0.25pt width\xcoord \Einheit depth0.25pt\hss}}
\def\Gitter(#1,#2){\unskip \xcoord0 \ycoord0 \leavevmode
  \LOOP\ifnum\ycoord<#2
    \loop\ifnum\xcoord<#1
      \raise\ycoord \Einheit\hbox to0pt{\hskip\xcoord 
\Einheit\Punkt\hss}%
      \advance\xcoord by1
    \repeat
    \xcoord0
    \advance\ycoord by1
  \REPEAT}
\def\Gitter(#1,#2)(#3,#4){\unskip \xcoord#3 \ycoord#4 \leavevmode
  \LOOP\ifnum\ycoord<#2
    \loop\ifnum\xcoord<#1
      \raise\ycoord \Einheit\hbox to0pt{\hskip\xcoord 
\Einheit\Punkt\hss}%
      \advance\xcoord by1
    \repeat
    \xcoord#3
    \advance\ycoord by1
  \REPEAT}
\def\Label#1#2(#3,#4){\unskip \xdim#3 \Einheit \ydim#4 \Einheit
  \def\lo{\advance\xdim by-.5 \Einheit \advance\ydim by.5 \Einheit}%
  \def\llo{\advance\xdim by-.25cm \advance\ydim by.5 \Einheit}%
  \def\loo{\advance\xdim by-.5 \Einheit \advance\ydim by.25cm}%
  \def\o{\advance\ydim by.25cm}%
  \def\ro{\advance\xdim by.5 \Einheit \advance\ydim by.5 \Einheit}%
  \def\rro{\advance\xdim by.25cm \advance\ydim by.5 \Einheit}%
  \def\roo{\advance\xdim by.5 \Einheit \advance\ydim by.25cm}%
  \def\l{\advance\xdim by-.30cm}%
  \def\r{\advance\xdim by.30cm}%
  \def\lu{\advance\xdim by-.5 \Einheit \advance\ydim by-.6 \Einheit}%
  \def\llu{\advance\xdim by-.25cm \advance\ydim by-.6 \Einheit}%
  \def\luu{\advance\xdim by-.5 \Einheit \advance\ydim by-.30cm}%
  \def\u{\advance\ydim by-.30cm}%
  \def\ru{\advance\xdim by.5 \Einheit \advance\ydim by-.6 \Einheit}%
  \def\rru{\advance\xdim by.25cm \advance\ydim by-.6 \Einheit}%
  \def\ruu{\advance\xdim by.5 \Einheit \advance\ydim by-.30cm}%
  #1\raise\ydim\hbox to0pt{\hskip\xdim
     \vbox to0pt{\vss\hbox to0pt{\hss$#2$\hss}\vss}\hss}%
}
\catcode`\@=12

\usepackage{amsmath,amsthm}
\usepackage{amssymb}
\usepackage{color}
\usepackage{xspace}
\usepackage[colorlinks=true,
linkcolor=green,
filecolor=brown,
citecolor=green]{hyperref}

\def\v{\vert}

\def\u{\ensuremath{\mathcal U}\xspace}
\def\rl{R-L minima\xspace}

\def\mbf#1{\mathchoice{\hbox{\boldmath $\displaystyle #1$}}
        {\hbox{\boldmath $\textstyle #1$}}
        {\hbox{\boldmath $\scriptstyle #1$}}
        {\hbox{\boldmath $\scriptscriptstyle #1$}}} % mathboldface

\begin{document}
\newtheorem{theorem}{Theorem}
\newtheorem{defn}[theorem]{Definition}
\newtheorem{lemma}[theorem]{Lemma}
\newtheorem{prop}[theorem]{Proposition}
\newtheorem{cor}[theorem]{Corollary}
%\vspace*{-5mm}
\begin{center}
{\Large
Pattern Avoidance in ``Flattened'' Partitions                         \\ 
}

\vspace{10mm}
DAVID CALLAN  \\
Department of Statistics  \\
\vspace*{-1mm}
University of Wisconsin-Madison  \\
\vspace*{-1mm}
1300 University Ave  \\
\vspace*{-1mm}
Madison, WI \ 53706-1532  \\
{\bf callan@stat.wisc.edu}  \\
\vspace{5mm}

February 15, 2008
\end{center}

\begin{abstract}
To flatten a set partition 
(with apologies to \emph{Mathematica$^{\textrm{\textregistered}}$}) means 
to form a permutation by erasing the dividers between its blocks. Of 
course, the result depends on how the blocks are listed. For the usual 
listing---increasing entries in each block and blocks arranged 
in increasing order of their first entries---we count the partitions 
of $[n]$ whose flattening avoids a single 3-letter pattern. Five 
counting sequences arise: a null sequence, the powers of 2, the Fibonacci numbers, the Catalan 
numbers, and the binomial transform of the Catalan numbers. 
\end{abstract}

\vspace{5mm}
\section{Introduction} 
There is an extensive literature on pattern avoidance in 
permutations. Klazar \cite{kla1,kla2,kla3} considered an analogous notion 
for set partitions and Sagan \cite{sagan06} introduced a second such 
notion based on restricted growth functions (see also 
\cite{goyt06,jelinek07}). Here we consider set partitions avoiding a 
permutation in the following sense. Suppose partitions $\Pi$ of 
$[n]=\{1,2,\ldots,n\}$ are written in some pre-specified standard form, 
say \emph{standard increasing} form: increasing entries in each block and blocks arranged 
in increasing order of their first entries. Then define Flatten($\Pi)$ to be 
the permutation of $[n]$ obtained by erasing the dividers between the 
blocks of $\Pi$. For example, $\Pi = $
136--279--4--58 is in standard increasing form and 
Flatten($\Pi)=136279458$. (The computer algebra system 
\emph{Mathematica$^{\textrm{\textregistered}}$} implements this 
operation with the command \texttt{Flatten}).
For a permutation $\pi$ on an initial segment of the positive integers (a 
pattern permutation) we say the partition $\Pi$ avoids $\pi$ or $\Pi$ 
is $\pi$-avoiding if the permutation
Flatten($\Pi)$ avoids $\pi$ (in the classical sense). We write $\Pi \vdash [n]$ if $\Pi$ is 
a partition of $[n]$. Set  $\u(n;\pi) =\{\Pi \vdash [n] :$ Flatten($\Pi)$ avoids $\pi$\}.

In \S 2, we fix standard increasing as the form for writing 
partitions of $[n]$ and count $\u(n;\pi)$
for all 3-letter pattern permutations $\pi$.

\vspace{5mm}
\section{Set partitions in standard increasing form} 
\subsection{123-avoiding}
This case is not very interesting; the counting sequence $(\v 
\u(n;123)\v)_{n\ge 1}$ is $(1,2,1,0,0,0,\ldots)$.

\subsection{132-avoiding}
A partition $\Pi$ of $[n]$ is in $\u(n;132)$ if and only if 
Flatten($\Pi$) is the identity permutation. This is because the first 
entry of Flatten($\Pi$) is always 1 and will be the `1' of a 132 
pattern unless Flatten($\Pi$) is an increasing sequence, that is, the 
identity permutation. So any subset of the $n-1$ 
spaces between $1,2,\ldots,n $ can serve as the dividers to form 
$\Pi$ and $\v \u(n;132) \v =2^{n-1}$.

\subsection{213-avoiding}
First, we claim a partition $\Pi$ of $[n]$ is in $\u(n;213)$ if and only if (i) the 
first block of $\Pi$ has the form $I \cup J$ with $I$ a nonempty 
initial segment of $[n]$ and $J$ a terminal segment of $[n]$ (possibly empty) disjoint 
from $I$, and (ii) the remaining blocks, when standardized, 
themselves form a 
$213$-avoiding partition. (To standardize means to replace smallest 
entry by 1, second smallest by 2, and so on.) 

Clearly, these two conditions are sufficient and condition (ii) is 
necessary. If condition (i) fails for $\Pi \in \u(n;213)$, let $a$ be 
the smallest element of $[n]$ not in the first block; $a$ is necessarily 
the first element of the second block. Because the condition fails 
there exist $b,c$ in $[n]$ with $c>b>a ,\ b$ in the first block and $c$ in 
a later block. Hence $c$ occurs after $a$ and $bac$ is a 213-pattern in 
Flatten($\Pi$), a contradiction. So condition (i) is necessary also. 

Now let $u(n) =\v\, \u(n;213)\,\v$ and set $u(n,k) =\v \, \{\Pi\in\u(n;213):\:$ 
first block of $\Pi$ has 
length $k\}\,\v$. Clearly, $u(n,n)=1$ and for $1\le k \le n$, the first 
block is determined by $I$ and there are $k$ choices for $I$, namely, 
$\big( [i] \big)_{i=1}^{k}$. Hence we have the system of equations
\begin{alignat*}{2}
    u(n,n) & = \quad 1 & &\quad \textrm{for $n\ge 1$}  \\
    u(n,k) & =  ku(n-k) & &\quad \textrm{for $1 \le k < n$}  \\
    u(n) & =  \sum_{k=1}^{n}u(n,k) &  &\quad \textrm{for $n\ge 1$} 
\end{alignat*}
with solution involving the Fibonacci numbers ($F_{-1}:=1,\, 
F_{0}=0,\,F_{1}=1$)
\begin{alignat*}{2}
    u(n,j) & = jF_{2n-2j-1} & &\quad \textrm{for $1 \le j < n$}  \\
    u(n) & =  F_{2n-1}. &  &
\end{alignat*}

\subsection{231-avoiding}
This case gives rise to the Catalan numbers via Touchard's identity 
\cite{shapiro76},
\begin{equation}
    C_{n}=\sum_{k\ge 0}\binom{n-1}{2k}2^{n-1-2k}C_{k}.
    \label{touchard1}
\end{equation}
For a permutation $p$ of $[n]$, a \emph{descent terminator} is an entry smaller than its 
immediate predecessor and, by convention, the first entry is also 
considered a descent terminator. A \emph{right-to-left \emph{(\emph{R-L})} minimum} of $p$ is an entry 
smaller than all the entries after it. Clearly, for a partition in 
standard increasing form and its associated permutation, $\{\textrm{descent 
terminators}\} \subseteq \{\textrm{block initiators}\} \subseteq \{\textrm{\rl}\}$.
For $\Pi\vdash [n]$, let $M(\Pi)$ denote the set of \rl of 
Flatten($\Pi$) that are not descent terminators, and set 
$\u(n,k\, ;231)=\{\Pi \in \u(n;231):\,\v M(\Pi) \v =k\}$. We claim $\v \,
\u(n,k\, ;231)\,\v =\binom{n-1}{k}2^{k}C_{\frac{n-1-k}{2}}$ where $C_{n}$ 
is the Catalan number and $C_{n}:=0$ when $n$ is not an integer. 
Touchard's identity (\ref{touchard1}) then implies $\v\, 
\u(n;231)\,\v=C_{n}$.

To establish the claim, it suffices to show
\begin{eqnarray}
    \v\, \u(n,0\, ;231)\, \v  & = &  C_{\frac{n-1}{2}} \qquad \textrm{for 
    $n\ge 1$, and}
    \label{A1}  \\[4pt]
   \v\, \u(n,k\, ;231)\, \v & = & \binom{n-1}{k}
    2^{k}\v\,\u(n-k,0\, ;231)\,\v \quad \textrm{for 
    $1 \le k < n$.} 
    \label{A2}
\end{eqnarray}
\quad To show (\ref{A1}), let $\Pi\in \u(n,0\, ;231)$. Then the \rl and descent terminators 
of $p:=$ Flatten($\Pi$) coincide. 
The last entry of $p$ is certainly a 
R-L minimum, hence a descent terminator, and so it must form a 
singleton block in $\Pi$.
Each non-last block has length $\le 2$ because if $(a,b,c,\ldots)$ is a 
block of length $\ge 3$, then $bcd$ is a 231-pattern where $d$ is the 
first entry of the next block: certainly $b<c$ and we also have $d<b$ 
because if $b<d$, then $b<$ \emph{all} entries that follow it. This 
would make $b$ a R-L minimum that was not a descent terminator, a 
contradiction.
On the other hand, each non-last block has length $\ge 2$ because a 
non-last singleton block would imply that the first entry of the next 
block was a R-L minimum that was not a descent terminator. Hence all 
but the last block have length 2 and so $n$ is odd, say $n=2r+1$, 
and $\Pi$ is of the form 
$a_{1}b_{1}$\,--\,$a_{2}b_{2}$\,--\,\ldots\,--\,$a_{r}b_{r}$\,--\,$a_{r+1}$. 

Clearly, 
$a_{1}=1$. Also, $a_{2}=2$ because otherwise, since $a_{2}$ is a R-L 
minimum, 2 would occur to the left of $a_{2}$ and this would force 
$b_{2}=2$. But then $a_{2}$ would be a R-L minimum that was not a descent terminator.
Next, we claim $a_{i+2}\le 2i+2$ for $1\le i \le r-1$. Suppose 
contrariwise that $a_{i+2}>2i+2$ for some $i$. Then none of 
$3,4,\ldots,2i+2$ can occur after $a_{i+2}$ because $a_{i+2}$ is a R-L 
minimum. This forces the first $i+1$ blocks to consist of the first 
$2i+2$ positive integers leaving $b_{i+1}$ a R-L minimum, which 
is not possible. Hence the sequence $(c_{i})_{i=1}^{r-1}$ with 
$c_{i}:=a_{i+2}-2$ satisfies 
\begin{equation}
    \label{item:t}
1 \le c_{1}<c_{2}<\ldots <c_{r-1}, \quad \textrm{and} \quad
c_{i}\le 2i \quad \textrm{for $1\le i \le r-1$}.
\end{equation}

We have exhibited a map from $\u(2r+1,0\, ;231)$ to sequences 
$(c_{i})_{i=1}^{r-1}$ satisfying (\ref{item:t}). This map is in fact 
a bijection and here is its inverse. Given such a sequence, for 
example with $r=9,\ (c_{i})=(1,2,4,5,7,12,13,15)$, we can immediately recover 
the $a_{i}$'s and must determine the $b_{i}$'s (blank squares in Fig. 1).

\vspace*{-5mm}

\[
\begin{array}{cccccccccc}
a_{1}\ b_{1} &   a_{2}\ b_{2} & \ldots & & & &  &\ldots  & a_{r}\ b_{r} & a_{r+1} \\
1\ \square & 2\ \square & 3\ \square & 4\ \square & 6\ \square & 7\ \square & 
9\ \square & 14\ \square & 15\ \square & 17
 \end{array}
\]
\centerline{Fig. 1}

\vspace*{1mm}

\noindent Fill in the blank squares using 
$B=[2r+1]\backslash(a_{i})_{i=1}^{r+1}$ from right to left as 
follows. Place the smallest element of $B$ that exceeds 
$a_{r+1}$ in the $b_{r}$ square and, in general, place the smallest 
not-yet-placed element of $B$ that exceeds $a_{i+1}$ in the $b_{i}$ 
square. The example has $B=\{5,8,10,11,12,13,16,18,19\}$, yielding 
$(b_{i})_{i=1}^{r}=(13,12,5,11,8,10,19,16,18)$.

There is a nice graphical way to visualize the result of this 
algorithmic procedure using Dyck paths.
Recall that the 
Catalan number $C_{r}$ counts sequences $(c_{i})_{i=1}^{r-1}$ satisfying 
(\ref{item:t}) \cite[Ex.\,6.19, item t]{ec2}.
Indeed, given a Dyck path of semilength $r$ let $c_{i}$ denote the 
number of steps preceding the $(i+1)$st upstep for $1\le i \le r-1$. 
This is a bijection from Dyck $r$-paths to the sequences 
$(c_{i})_{i=1}^{r-1}$ satisfying (\ref{item:t}). So, sketch the Dyck 
path corresponding to the sequence $(c_{i})_{i=1}^{r-1}$, prepend an 
upstep, and number all $2r+1$ steps in order left to right, as in 
Fig.2 for our running example.

\vspace*{-5mm}

\Einheit=0.6cm
\[
\Pfad(-10,0),3333433434444334344\endPfad
\SPfad(-9,1),111111111111111111\endSPfad
\DuennPunkt(-10,0)
\DuennPunkt(-9,1)
\DuennPunkt(-8,2)
\DuennPunkt(-7,3)
\DuennPunkt(-6,4)
\DuennPunkt(-5,3)
\DuennPunkt(-4,4)
\DuennPunkt(-3,5)
\DuennPunkt(-2,4)
\DuennPunkt(-1,5)
\DuennPunkt(0,4)
\DuennPunkt(1,3)
\DuennPunkt(2,2)
\DuennPunkt(3,1)
\DuennPunkt(4,2)
\DuennPunkt(5,3)
\DuennPunkt(6,2)
\DuennPunkt(7,3)
\DuennPunkt(8,2)
\DuennPunkt(9,1)
\Label\o{ \textrm{{\scriptsize 1}}}(-9.8,0.2)
\Label\o{ \textrm{{\scriptsize 2}}}(-8.8,1.2)
\Label\o{ \textrm{{\scriptsize 3}}}(-7.8,2.2)
\Label\o{ \textrm{{\scriptsize 4}}}(-6.8,3.2)
\Label\o{ \textrm{{\scriptsize 5}}}(-5.7,2.9)
\Label\o{ \textrm{{\scriptsize 6}}}(-4.8,3.2)
\Label\o{ \textrm{{\scriptsize 7}}}(-3.8,4.2)
\Label\o{ \textrm{{\scriptsize 8}}}(-2.7,3.9)
\Label\o{ \textrm{{\scriptsize 9}}}(-1.8,4.2)
\Label\o{ \textrm{{\scriptsize 10}}}(-0.8,3.9)
\Label\o{ \textrm{{\scriptsize 11}}}(0.2,2.9)
\Label\o{ \textrm{{\scriptsize 12}}}(1.2,1.9)
\Label\o{ \textrm{{\scriptsize 13}}}(2.2,0.9)
\Label\o{ \textrm{{\scriptsize 14}}}(3.1,1.2)
\Label\o{ \textrm{{\scriptsize 15}}}(4.1,2.2)
\Label\o{ \textrm{{\scriptsize 16}}}(5.2,1.9)
\Label\o{ \textrm{{\scriptsize 17}}}(6.2,2.2)
\Label\o{ \textrm{{\scriptsize 18}}}(7.2,1.9)
\Label\o{ \textrm{{\scriptsize 19}}}(8.2,0.9)
\]

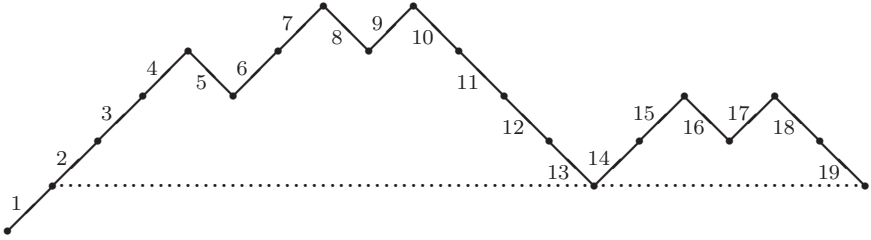
\centerline{Fig. 2}

\vspace*{1mm}

\noindent Every upstep in a Dyck path has a matching downstep: the first one 
encountered directly east from the upstep or, more precisely, the 
terminal downstep of the shortest Dyck subpath starting at the upstep. 
The $a_{i}$'s are evident in the augmented Dyck path as the labels on 
the upsteps, and the $b_{i}$'s are also discernible: $b_{i}$ is is 
the label on the matching downstep for the \emph{next} upstep after $a_{i},\ 
1\le i \le r$. It is now clear that the $a_{i}$'s are increasing and that
$a_{i}<b_{i}>a_{i+1}$ for $1\le i \le r$; hence $(a_{i})_{i=1}^{r+1}$ is both 
the set of \rl and the set of descent terminators in Flatten($\Pi$) 
and so $\Pi \in \u(2r+1,0\, ;231)$. It is also easy to verify that $\Pi$ 
is 231-avoiding. Indeed, since all entries following $a_{i}$ 
are $> a_{i}$, the first two entries of a putative 231 pattern would 
have to be $b$'s, say $b_{i}<b_{j}$ with $i<j$, and $a_{j}$ would be 
the last upstep preceding $b_{i}$ (or else $b_{j}$ would be $< 
b_{i}$). Hence, for all $k>j$, upstep $a_{k}$ occurs after $b_{i}$ and 
so $b_{k}>a_{k}>b_{i}$ for $k>j$. Since $b_{i}$ is the `2' of the 
231 pattern and we have just seen that all later entries are larger 
than $b_{i}$, no entry after $b_{j}$ can serve as the `1' of the 
pattern. We conclude that the partition 
$a_{1}b_{1}$\,--\,$a_{2}b_{2}$\,--\,\ldots\,--\,$a_{r}b_{r}$\,--\,$a_{r+1}$ is in $\u(n,0\, ;231)$ as required.

To prove (\ref{A2}), consider $\Pi \in \u(n,k\, ;231)$. Let $K$ denote the 
set of \rl that are not descent terminators in Flatten($\Pi$). Thus 
$\v\,K\,\v = k$  and $K \subseteq [2,n]$. Let $L$ denote the set of 
elements in $K$ that initiate a block in $\Pi$. Thus $L\subseteq K$. 
Let $\Pi_{0}$ denote the partition obtained from $\Pi$ by deleting 
each element $i$ of $K$ from its block and, if $i$ is also in $L$,
concatenating this block with the currently preceding block. Then $\Pi_{0} \in 
\u(n-k,0\,;231)$. 
For example, $\Pi = $ 1\,--\,24\,--\,37\,--\,568 yields 
$K=\{2,6,8\},\ L=\{2\}$, and $\Pi_{0}= $ standardize(14\,--\,37\,--\,5) = 
13\,--\,25\,--\,4. An example where three consecutive blocks are 
concatenated to form $\Pi_{0}$ is  $\Pi = $ 1\,--\,2\,--\,35\,--\,4 with 
$K=\{2,3\},\ L=\{3\}$, and $\Pi_{0}= $ standardize(15\,--\,4) = 
13\,--\,2. We claim the map\ \ $\u(n,k\, ;231) \longrightarrow 
(K,L,\Pi_{0})$ is a bijection to all triples $(K,L,\Pi_{0})$ with $K$ 
a $k$-element subset of $[2,n]$, $L$ an arbitrary subset of $K$, 
and $\Pi_{0}$ a partition in $\u(n-k,0\,;231)$, and (\ref{A2}) then follows 
from (\ref{A1}). To establish the claim, suppose given such a triple 
$(K,L,\Pi_{0})$, and build up $\Pi$ as follows from $\Pi_{0}$. For 
each $a\in K$ in turn from smallest to largest, locate the last block 
in the current partition whose first entry is $<a$; then, to get the 
next partition, after adding 1 to each entry $\ge a$ insert $a$ into 
the located block at the appropriate position to ensure an increasing 
block. The end result will be a partition of $[n]$ in which the 
descent terminators are the block initiators and no element of $K$ is 
a block initiator. Finally, for each element of $K$ that is in the 
subset $L$, place a divider just before that element so that it 
initiates a block. This procedure yields $\Pi$ and shows the map is 
invertible.

\subsection{312-avoiding}
We claim a partition $\Pi$ of $[n]$ is in $\u(n;312)$ if and only if (i) the 
first block of $\Pi$ is all of $[n]$ or has the 
form $I\backslash \{a\}$ where $I$ is an 
initial segment of $[n]$ of length $\ge 2$ and $a\ge 2$ is in $I$, and (ii) the remaining blocks, when standardized, 
themselves form a $312$-avoiding partition. 

The conditions are 
sufficient because if they hold and a 312 pattern involved the first 
block, then only the `3' could occur in the first block leaving 
the `1' and `2' to occur in later blocks. This however is 
impossible because at most one letter smaller than the `3' is missing 
from the first block. So we merely need to show that condition (i) is 
necessary. Suppose then that  condition (i) is not met. Let $c$ denote 
the largest entry in the first block and $a$ the smallest letter 
missing from the first block. Then by supposition there is a letter 
$b$ missing from the first block with $a<b<c$. Since $a$ must be the 
first entry of the second block, $b$ occurs after $a$ and $cab$ is a 
312 pattern in Flatten($\Pi$), a contradiction.

Now, if the first block has length $k<n$, there are exactly $k$ choices 
for $a$, namely, $2,3,\ldots,k+1$. This observation leads to the very same 
recurrence relation as in the 213-avoiding case, and another 
Fibonacci counting sequence: $\v\, \u(n;312)\,\v =F_{2n-1}$.

\subsection{321-avoiding} 
This case is counted by the binomial transform of the Catalan 
numbers: $\v\, \u(n+1;321)\,\v= \sum_{k=0}^{n}\binom{n}{k}C_{k}$. Our
proof is quite similar to that of the 231-avoiding case but with 
Touchard's identity replaced by the following one involving the Riordan 
numbers $R_{n}$, 
\begin{equation} 
    \sum_{k=0}^{n}\binom{n}{k}2^{k}R_{n-k}=\sum_{k=0}^{n}\binom{n}{k}C_{k},
    \label{riordan}
\end{equation}
where $R_{n}:=\sum_{j=0}^{n}(-1)^{n-j}\binom{n}{j}C_{j}$.
%(A more direct proof would be desirable.)
The identity (\ref{riordan}) is easily proved by reversing the order 
of summation after substituting for $R_{n-k}$.

The Riordan number $R_{n}$ 
(\htmladdnormallink{A005043}{http://www.research.att.com:80/cgi-bin/access.cgi/as/njas/sequences/eisA.cgi?Anum=A005043}
in OEIS) is well known to count, among other 
things, Dyck $n$-paths with no short descents. (A `descent' is a 
maximal sequence of contiguous downsteps and `short' means of length 
1.) 
Mimicking Section 2.4, define $\u(n,k\, ;321)=\{\Pi \in \u(n;321):\,\v M(\Pi) \v =k\}$. 
We claim $\v \,
\u(n,k\, ;321)\,\v =\binom{n-1}{k}2^{k}R_{n-1-k}$ for $0\le k \le 
n-1$, and the identity (\ref{riordan}) then implies $\v\, 
\u(n;321)\,\v=\sum_{k=0}^{n-1}\binom{n-1}{k}C_{k}$.

To establish the claim, it suffices to show
\begin{eqnarray}
    \v\,\u(n,0\,;321)\,\v & = & R_{n-1} 
    \label{B1}  \\[4pt]
   \v\,\u(n,k\,;321)\,\v & = & \binom{n-1}{k}2^{k} \v\,\u(n-k,0\,;321)\,\v 
   \quad \textrm{for $1\le k \le n-1$}
    \label{B2}
\end{eqnarray}
To prove assertion (\ref{B1}) there is a bijection (essentially due to 
Krattenthaler \cite{krattenthaler}) from $\u(n,0\,;321)$ to Dyck 
$(n-1)$-paths with no short descents, illustrated with $n=9$:

%\vspace*{-2mm}

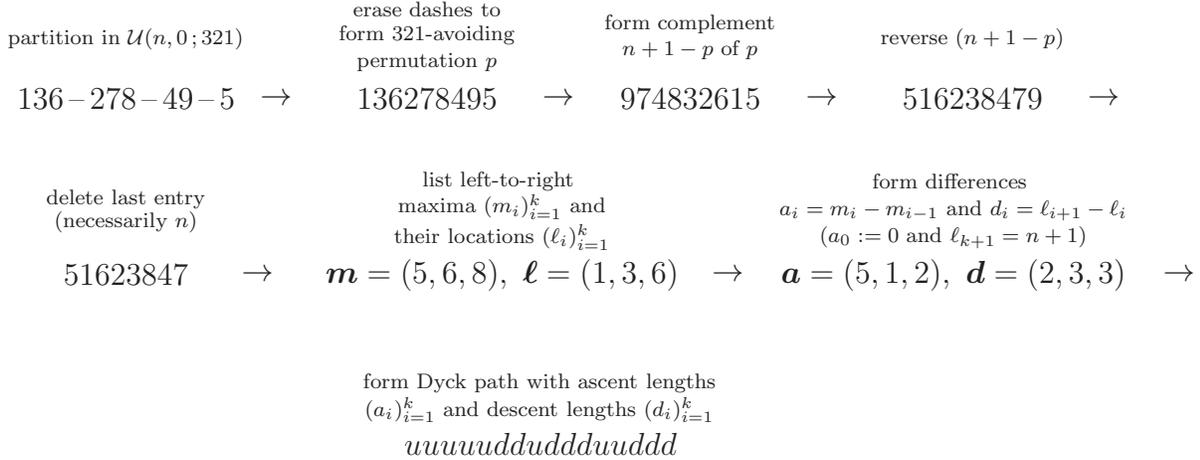
\Einheit=0.5cm
\[
\Label\o{ \textrm{{\scriptsize partition in $\mathcal{U}(n,0\,;321)$}}}(-13,3.3)
\Label\o{ \textrm{136\,--\,278\,--\,49\,--\,5}}(-13,1.7)
\Label\o{ \textrm{{\scriptsize erase dashes to}}}(-5,4.1)
\Label\o{ \textrm{{\scriptsize form 321-avoiding}}}(-5,3.4)
\Label\o{ \textrm{{\scriptsize permutation $p$}}}(-5,2.7)
\Label\o{ \textrm{136278495}}(-5,1.7)
\Label\o{ \textrm{{\scriptsize form complement}}}(2,3.7)
\Label\o{ \textrm{{\scriptsize $n+1-p$ of $p$}}}(2,3.0)
\Label\o{ \textrm{974832615}}(2,1.7)
\Label\o{ \textrm{{\scriptsize reverse ($n+1-p$)}}}(9.5,3.3)
\Label\o{ \textrm{516238479}}(9.5,1.7)
\Label\o{ \rightarrow}(-9,1.7)
\Label\o{ \rightarrow}(-1.5,1.7)
\Label\o{ \rightarrow}(5.5,1.7)
\Label\o{ \rightarrow}(13,1.7)
\]

\vspace*{-10mm}

\[
\Label\o{ \textrm{{\scriptsize delete last entry}}}(-13,3.7)
\Label\o{ \textrm{{\scriptsize (necessarily $n$)}}}(-13,3.1)
\Label\o{ \textrm{51623847}}(-13,1.7)
\Label\o{ \textrm{{\scriptsize list left-to-right }}}(-3,4.2)
\Label\o{ \textrm{{\scriptsize maxima $(m_{i})_{i=1}^{k}$ and}}}(-3,3.5)
\Label\o{ \textrm{{\scriptsize their locations 
$(\ell_{i})_{i=1}^{k}$}}}(-3,2.7)
\Label\o{ \textrm{$\mbf{m}=(5,6,8),\ \mbf{\ell}=(1,3,6)$}}(-3,1.7)
\Label\o{ \textrm{{\scriptsize form differences }}}(9,4.2)
\Label\o{ \textrm{{\scriptsize $a_{i}=m_{i}-m_{i-1}$ and 
$d_{i}=\ell_{i+1}-\ell_{i}$}}}(9,3.4)
\Label\o{ \textrm{{\scriptsize ($a_{0}:=0$ and $\ell_{k+1}=n+1$)}}}(9,2.7)
\Label\o{ \textrm{$\mbf{a}=(5,1,2),\ \mbf{d}=(2,3,3)$}}(9,1.7)
\Label\o{ \rightarrow}(-9.5,1.7)
\Label\o{ \rightarrow}(3,1.7)
\Label\o{ \rightarrow}(15,1.7)
\Label\o{ \textrm{{\scriptsize form Dyck path with ascent 
lengths}}}(-2,-1.2)
\Label\o{ \textrm{{\scriptsize $(a_{i})_{i=1}^{k}$ and descent lengths 
$(d_{i})_{i=1}^{k}$}}}(-2,-1.9)
\Label\o{ \textrm{$uuuuudduddduuddd$}}(-2,-2.8)
\]

\vspace*{5mm}

\centerline{{\small Bijection $\u(n,0\,;321) \longrightarrow$ to Dyck 
$(n-1)$-paths with no short descents}}

\vspace*{5mm}

The proof of assertion (\ref{B2}) uses the same bijection $\Pi 
\rightarrow (K,L,\Pi_{0})$ as in the proof of 
(\ref{A2}) and is omitted.

It would be interesting to investigate permutation-avoidance for other canonical 
representations of a set partition where less familiar counting 
sequences seem to arise.

\end{document}

\htmladdnormallink{\emph{Electronic J. Combin.}}{http://www.combinatorics.org/}

\begin{alignat*}{2}
    (i)\qquad & \v\, \u(n,0\, ;231)\, \v = C_{\frac{n-1}{2}} & &\quad \textrm{for 
    $n\ge 1$, and}  \\[5pt] 
    (ii)\qquad  &\v\, \u(n,k\, ;231)\, \v =\binom{n-1}{k}
    2^{k}\v\,\u(n-k,0\, ;231)\,\v & &\quad \textrm{for 
    $1 \le k < n$.} 
\end{alignat*}